\documentclass[graybox]{SNmult}
\usepackage{type1cm}    
\usepackage{makeidx}         
\usepackage{graphicx}       
\usepackage{multicol}      
\usepackage[bottom]{footmisc}
\usepackage{newtxtext}
\usepackage[varvw]{newtxmath}
\makeindex           

\begin{document}

	\title*{Steiner porism and boundary value problems}

	\author{Giorgi Khimshiashvili} \institute{Ilia State University, Ilia Chavchavadze Ave. N45, 
		Tbilisi 0179, Georgia; \email{giorgi.khimshiashvili@iliauni.edu.ge}}
	
	\maketitle
	\abstract{The aim of this paper is to describe several interpretations of Steiner chains of circles in terms of certain boundary value problems and give a short review of related results. Our approach and discussion are motivated by the recent results on the geometry of Steiner porism based on an apparently novel analytic parametrization of poristic circles obtained by A.Diakvnishvili and the present author. We also mention a few feasible problems and generalizations suggested by the constructions and considerations presented in the paper.}
	
	
	    
	\section{Introduction}
	\label{sec:1}
	
	We discuss certain connections between the classical Steiner porism and certain boundary value problems of Riemann-Hilbert type. Our approach is based on an analytic representation of poristic Steiner chains developed recently in \cite{bisa}, \cite{bidi}, \cite{dikh}, which yields a loop in the so-called Brylinski loop space of three-dimensional Euclidean space $\mathbb{R}^3$ (see \cite{bry} for the definition and basic properties of Brylinski loop spaces). This in turn yields interpretations of Steiner chains in terms of the nonlinear Riemann-Hilbert problem (RHP) in the sense of \cite{weg}, in a general setting of Riemann-Hilbert problem in almost complex manifolds developed in \cite{khi1}, \cite{khi2}, and in terms of boundary value problems for circle packings considered as discrete analytic functions \cite{ste}, \cite{glwe}. We give a short review of these topics which basically consists of corollaries of some related results presented in \cite{dikh}, \cite{glwe}, \cite{werokr} supplied with short comments on possible further developments. In doing so we do not pursue generality and freely use the necessary concepts and results from the papers mentioned above. 
	
	The central object of our discussion are the Steiner chains of circles in Euclidean plane. So let us recall that the term {\it Steiner $n$-chain} refers to a sequence of circles $\delta_i, i=1, \ldots, n,$ in Euclidean plane such that the circles with adjacent indices $mod\, n$ are externally tangent, each of $\delta_i$ is externally tangent to a fixed circle $\gamma$ and internally tangent to another fixed circle $\Gamma$ containing $\gamma$ in its interior disk. The nested circles $\gamma$ and $\Gamma$ are called the {\it parent circles} of the Steiner chain considered \cite{ber}. Steiner chains is a classical topic discussed in many papers, in particular, in the context of {\it Steiner porism} (see, e.g., \cite{ped}, \cite{scta}, \cite{bidi}) which plays a central role in our discussion. We will also deal with a similar but a bit more general concept of {\it circle garland} in the definition of which is not required the existence of parent circles. 
	
	The boundary value problems considered in the sequel arise in the contexts of nonlinear Riemann-Hilbert problems in the sense of \cite{weg}, holomorphic curves in Brylinski loop spaces \cite{khi2}, \cite{khsi}, \cite{khwo}, and circle packings \cite{ste}. The necessary background on these topics is given below. An essential part of our discussion relies on recent results on Steiner $n$-chains obtained in \cite{mazy}, \cite{bisa}, \cite{bidi}.

	\section{Preliminaries on Steiner chains and their invariants}
	\label{sec:2}
	
	The aforementioned Steiner porism states that there exists a one-dimensional family of Steiner $n$-chains having the same parent circles \cite{scta}. By a way of analogy with the Poncelet porism \cite{ber} a collection of Steiner $n$-chains with the fixed parent circles will be called {\it poristic Steiner $n$-chains}. 
	
	Let $n\geq 3$ be a natural number and $(R, r, d)$ with $R>r$ be a triple of positive numbers satisfying the relation 	
	\begin{equation} \label{pedoe-n}
		d^2 =(R-r)^2-4qRr,
	\end{equation}
	where $q=\tan^2\frac{\pi}{n}.$ It is known that such a triple defines a poristic system of Steiner $n$-chains       with the parent circles having radii $R, r$ and the distance between their centers equal to $d$ (see, e.g., \cite{ped}). We will say that this is a pair of parent circles of order $n$ with the gauge $(R,r,d)$. The number $d$ will be called the {\it gap} of parent circles. If the gap $d$ is equal to the radius $r$ of the inner parent circle of Steiner porism such a pair of parent circles  will be called {\it exact}. Clearly, this is equivalent to the relation $R = (4q+2)r$ with $q=\tan^2\frac{\pi}{n}.$ The invariants of the corresponding poristic Steiner chains only depend on $r$, which simplifies many issues concerned with Steiner chains including those discussed in the sequel.  
	
	Denote by $a = \frac{1}{r}, A=-\frac{1}{R}$ the so-called {\it bends} (signed curvatures) of the parent circles. Following \cite{scta}, for each Steiner $n$-chain $Z$ and natural number $k$, we define the $k$-th {\it $b$-moment} $I_k^{(n)}(Z)$ of chain $Z$ as
	\begin{equation} \label{bendingmoments1}
		I_k^{(n)}(Z) = \sum_{j=1}^n b_j^k,
	\end{equation}
	where $b_j$ are the bends (curvatures) of the circles in $Z.$ As was shown in \cite{scta}, the first $n-1$ $b$-moments $I_k^{(n)}(Z)$ are invariant in a poristic family of Steiner $n$-chains with the given parent circles, so we denote them simply by $I_k^{(n)} = I_k^{(n)}(R,r,d)$ and call $I_k^{(n)}$ the {\it invariant $b$-moments} of poristic Steiner chains. Clearly, the invariant $b$-moments can be expressed through radii $R, r$ of the parent circles or, equivalently, through their signed curvature $a, A$. For further reference, we give such formulas for $n=3$ and $n=4$.
	
	For $n=3$, explicit formulas for the first two $b$-moments $I_1, I_2$ follow directly from the Descartes circle theorem (see, e.g., \cite{scta}) and look as follows:
	\begin{equation} \label{moment31}
		I_1^{(3)} = b_1 + b_2 + b_3 = \frac{A+a}{2} = \frac{R-r}{2Rr},
	\end{equation}
	\begin{equation} \label{moment32}
		I_2^{(3)} = b_1^2 + b_2^2 + b_3^2 = \frac{A^2 + 6Aa + a^2}{8} = \\
		\frac{R^2 - 6Rr + r^2}{8R^2r^2}.\\		.
	\end{equation}
	
	It seems worthy of noting here that, as was shown in \cite{dikh}, the formula for $I_2^{(3)}$ given in page 239 of \cite{scta} is not correct and should be replaced by (\ref{moment32}). 
	
	Obviously definition of $b$-moments is applicable for circle garlands. Moreover, for a $n$-circle garland $G$ and non-negative integers $k, m$, the so-called {\it $z$-moment} $J_{k,m}(G)$ has been defined in \cite{lamawi} $G$ by the formula:
	\begin{equation} \label{complexmoments}
		J_{k,m}(G) = \sum b_i^k\cdot z_i^m,
	\end{equation}
	where $b_i$ are the bends as above and $z_i$ are the centers of the circles in the given $n$-garland $G$, considered as complex numbers by complexifying the reference plane. Clearly, $J_{k,0} = I_k$ so the $z$-moments are generalizations of $b$-moments. The collection of $z$-moments contains important information on a circle garland, in particular, it enables one to restore the shape of garland up to a rigid motion of the plane. This also refers to some issues concerned with the boundary values discussed below. An important result of \cite{lamawi} states that $z$-moments $J_{k,m}$ with $0\leq m \leq k \leq n-1$ are invariants of poristic Steiner $n$-chains. 
	
	We call such quantities the {\it invariant $z$-moments} of porism and give formulas for those invariants for $n=4$. The total amount of all invariant $z$-moments is $(n-1)(n+2)/2$ including $n-1$ $b$-moments $J_{k,0} = I_k$ and their exist algebraic relations between these invariants. Clearly, the quantities $J_{k,m}(G)$ are in general complex numbers. However the invariant $z$-moments have a hidden real structure, which is apparent from the formulas given below. Namely, all quantities $J_{k,m}$ with $0\leq m \leq k \leq n-1$ become simultaneously real for an appropriate embedding of the parent circles in the complex plane \cite{dikh}. The latter observation may be new since to the best of our knowledge this fact was never mentioned in the literature. We emphasize that this property only refers to invariant $z$-moments. We give now some of the explicit formulas for the invariant $z$-moments obtained in \cite{dikh}, to which we refer in Section \ref{sec:3} and Section \ref{sec:4}. 
	
	\begin{proposition} (\cite{dikh})
		For a given pair of parent circles of order four with the gauge $(R,r,d)$, we have
		
		\begin{equation} \label{mom11}
			J^{(4)}_{1,1} = \frac{2d}{r}=\frac{2\sqrt{A^{2}+6Aa+a^{2}}}{A},
		\end{equation}
		
		\begin{equation} \label{mom21}
			J^{(4)}_{2,1} = \frac{d\,(3R - 5r)}{2 R r}=\frac{3a+5A}{2A}\,\sqrt{A^{2}+6Aa+a^{2}},
		\end{equation}
		
		\begin{equation} \label{mom31}
			J^{(4)}_{3,1} = \frac{(5R - 3r)(R - 3r)\,d}{4 R^{2} r^{3}}=\frac{(5a+3A)(a+3A)}{4A}\,\sqrt{A^{2}+6Aa+a^{2}},
		\end{equation}
		
		\begin{equation} \label{mom22}
			J^{(4)}_{2,2} = \frac{3 d^{2}}{2r}=\frac{3\sqrt{A^{2}+6Aa+a^{2}}}{4A}.
		\end{equation}
	\end{proposition}
	
	\section{Steiner chains and Riemann-Hilbert problems}
	\label{sec:3}
	
	We pass to a connection between Steiner chains and nonlinear RHPs based on the observation that a Steiner chain $Z$ defines a loop $L(Z)$ in the space ${\mathbb L}_3 = L(\mathbb{R}^3)$ of smooth loops in $\mathbb{R}^3$. This observation, in turn, is based on the following two formulas given in \cite{dikh}. In spite of their simplicity these formulas are possibly new since we were not able them or equivalent statements in the existing literature. 
	Let us consider a system of poristic Steiner chains generated by parent circles with the gauge $(R, r, d)$. 
	
	\begin{proposition} \label{poristicradius} (\cite{dikh})
		The radius $r(t)$ of a poristic circle with the center having a given polar angle $t$ in the canonical coordinate system is given by:	
		\begin{equation} \label{keyformula1}	
			r(t) = \frac{R^2 - 2dR\cos t + d^2 - r^2}{2(R + r - d\cos t)}.
		\end{equation}
	\end{proposition}
	
	A short proof of this statement can be found in \cite{dikh}. The following result is a kind of converse statement and follows by resolving the equation (\ref{keyformula1}) with respect to the polar angle $t$ assuming that radius $r(t)$ is known. To avoid any misunderstanding, the given value of $r(t)$ is denoted below by $\rho$.
	
	\begin{proposition} \label{polarangle}
		For any $\rho \in [r_*, r^*],$ there exist exactly two poristic circles with the radii equal to $\rho$ and their polar angles $t_{\pm}(\rho)$ are given by  
		\begin{equation} \label{keyformula2}	
			t_{\pm}(\rho) = \pm \cos^{-1} \Big(\frac{R^2 - r^2 + d^2 - 2(R+r)\rho}{2d(R-\rho)}\Big).
		\end{equation}
	\end{proposition}
	It follows that, for any $t\in [0,2\pi]$, the poristic circle centered on a ray with polar angle $t$ coincides with the circle $\Gamma(z(t), r(t))$ of radius $r(t)$ with the center $z(t) = (R - r(t)) e^{{\imath}t}.$  In this way
	poristic circles yield a parameterized loop in the loop space of the reference plane. Finally, we define by $L(R,r)(t)$ a loop in ${\mathbb L_3}$ obtained as the result of rotating circle $\Gamma(z(t), r(t))$ by $\pi/2$ about the polar ray of its center $z(t)$ so that $L(R,r)(t)$ lies into a plane perpendicular to the reference plane of the porism. The union of all loops $L(R,r)(t)$ defines a smooth parametrized surface in $\mathbb R^3$ which we denote by $M(R,r)$ and call a {\it poristic target manifold}. Using the parameterization of $L(R,r)(t)$ given by formulas $(\ref{keyformula1}, \ref{keyformula2})$ it is easy to verify that $M(R,r)$ is smoothly diffeomorphic to a two-dimensional torus and belongs to a class of admissible target manifolds defined in \cite{weg}. Moreover, by its construction $M(R,r)$ has circular target curves so the results of \cite{glwe} are applicable in this case. Now, for $d=0$ the poristic target manifold $M(R,r)$ is isometric to a usual round torus. Hence it is foliated by the images of solutions to the corresponding RHP. This granted, the same argument as in the proof of Lemma 4.7 in \cite{glwe} yields the following conclusion. 
	
	\begin{proposition} \label{poristicRHP}
		For any $n\geq 3,$ there exists a positive number $\epsilon(n)$ such that, for any pair of parent circles of order $n$ with the gap $d$ smaller than $\epsilon(n)$, a nonlinear RHP with the poristic target manifold $M(R,r)$ is regularly solvable in the sense of \cite{weg}. 
	\end{proposition}
	
	The construction presented above suggests several further developments along the lines of \cite{glwe}. In particular, one may wish to investigate connections of the poristic target $M(R,r)$ with the Blaschke products and Nevanlinna parametrizations described in loc. cit. For $n=4$, the target manifolds $M(R,r)$ arising from an exact pair of parent circles appeared regular. So one may try to use the homotopy continuation method to show that all poristic RHP with the gap $d\leq r$ are regularly solvable. It seems also plausible to obtain exact analytic formulas for the solutions of RHP with poristic target manifold. We do not dwell on these developments for the reason of space and proceed by briefly discussing two other interpretations of Steiner chains in terms of boundary value problems mentioned in Section \ref{sec:1}. 
	
	Notice first that a nonlinear Riemann-Hilbert problem in the sense of \cite{weg} is a particular case of the problem of analytic discs in almost complex manifold considered in \cite{gro} and \cite{khi1}. Recall that, for a three-dimensional Riemannian manifold $M$, J.-L.Brylinski introduced and studied a natural almost complex structure $J_M$ on the Frechet manifold $L(M)$ of $C^{\infty}$-loops in $M$ \cite{bry}. For brevity, the arising almost complex Frechet manifold $(L(M), J_M)$ is denoted by $B(M)$ and called the {\it Brylinski loop space} $B(M)$ of $M$. M.Gromov's definition of pseudo-holomorphic curve in an almost complex manifold suggested in a celebrated paper \cite{gro} yields, in particular, a concept of pseudo-holomorphic disc in Brylinski loop space. These concepts have been used in \cite{khi1}, \cite{khi2} to formulate several boundary value problems (BVP) for pseudo-holomorphic functions with values in Brylinski loop spaces which included several generalizations of the classical RHP. Later on, certain concrete problems of such type were studied in more detail in \cite{khi3}, \cite{khsi}, \cite{khwo}. In particular, several constructions have been given in \cite{khi3}, \cite{khsi}, \cite{khwo} which yielded explicit solutions to such problems in several interesting cases. 
	
	It is now easy to realize that the construction of loop $L(R,r)$ given above enables one to consider a Riemann-Hilbert problem in the sense of \cite{khi2} with values in Brylinski loop space $\mathbb B_3$ of $\mathbb R^3.$ For small values of gap $d$, the solvability of this problem can be proved using topological methods described in \cite{khi3}. Interesting and feasible as it is, this line of development is beyond the scope of this paper and so it is delayed for future publications. 
	
	We conclude this section by a short discussion of a connection between Steiner chains and the boundary value problems for circle packings considered in \cite{ste}, \cite{werokr}. Recall that the boundary condition of {\it Beurling's boundary value problem for circle packings} (BPCP) studied in \cite{werokr} is given by $n$-tuple $R_b = \{r_i\}$ consisting of positive numbers interpreted as the radii of a circle garland $G$ in the reference plane. The solvability of a wide class of such BPCPs follows from the so-called Koebe-Andreev-Thurston theorem \cite{ste}. Moreover, the results of \cite{werokr} imply that rather detailed information on solutions of such a BPCP can be obtained in the case where the collection $R_b$ is given by the radii of a Steiner $n$-chain. Furthermore, an effective method of computing solutions of BPCP with boundary condition given by a Steiner $n$-chain became available in the framework of an algebraic approach to circle packings presented in a recent preprint \cite{mazy2}. In such a situation, a general paradigm of combinatorial geometry suggests that many properties of solutions to BPCP should be expressible through the $z$-moments of the boundary Steiner chain, which seems to be an interesting and feasible perspective. For $n=4$, this problem can be attacked using the formulas for the invariant $z$-moments given in Section \ref{sec:2}. 
	
	Thus one may expect that many properties of BPCP with a given boundary condition $R_b$ would become tractable if we were able to establish that $R_b$ coincides with a collection of radii of a Steiner. The latter issue, known as the {\it feasibility problem for Steiner chains} (cf. \cite{kir}), was effectively solved for $n=4$ in \cite{bidi}. With a view to further generalizations, we conclude this paper by presenting a simpler solution based on the formulas for the invariant $z$-moments given in Section \ref{sec:2}. 
	
	\section{Complex invariants and feasibility problems for Steiner chains}
	\label{sec:4}
	
	The aforementioned feasibility problems for Steiner chains arise in the process of searching for analogous results valid for arbitrary circle garlands. One of the natural approaches to solving such problems uses the results on the complex moments of Steiner chains given in Section \ref{sec:2}.
	
	The feasibility criterion is given in \cite{bidi} in terms of complex moments of circle garland $G$ defined above. 
	For $n=4$, the feasibility problem considered reads as follows. Given an ordered quadruple of positive numbers $(r_1, r_2, r_3, r_4)$ find out if there exists a Steiner $4$-chain with such radii of its circles in the given order. To the best of our knowledge, this natural problem was considered for the first time in \cite{kir}, where the author gave its solution in a concrete case but did not discuss the general case. A general algorithmic solution to this problem was given in \cite{bidi}. We reproduce it here in a simplified version based on the formulas given in Section \ref{sec:2}. 
	
	Recall that the algorithm given in \cite{bidi} involved several steps. First, one computes the first three moments $(I_1, I_2, I_3)$ of the quadruple $b_i = 1/r_i$ called the {\it actual $b$-moments}. Then assuming the existence of a sought Steiner $4$-chain one uses the formulas for the first two $b$-moments given in \cite{bidi} to find the "virtual radii" $(\tilde R, \tilde r)$ of the "virtual parent circles". From the mentioned formulas follows that it is sufficient to solve the system
	$${2(a+A)= I_1, 3A^2 + 10Aa + 3a^2 = 2I_2},$$
	which reduces to the following quadratic equation
	\begin{equation} \label{parentradii}
		16a^2 - 8I_1a + (8I_2 - 3I_1^2) = 0.
	\end{equation} 
	So the solutions of the above system can be expressed by explicit formulas, which yields the virtual parent radii $(\tilde R, \tilde r)$. If $\tilde R^2 - 6 \tilde R \tilde r + \tilde r^2 > 0$ then the pair $(\tilde R, \tilde r)$ is a feasible candidate for the radii of parent circles having the distance between their centers equal to $\tilde d = \sqrt{\tilde R^2 - 6 \tilde R \tilde r + \tilde r^2 }.$ Next, one has to verify that the given values of radii $r_i$ belong to the segment $[\tilde r_*, \tilde r^*]$. The final step of the algorithm given in \cite{bidi} required computing the "virtual bending moments" and comparing them with the the actual $b$-moments. In view of our results this last step can be substituted by verifying that the actual $z$-moments of bends satisfy the relations (\ref{mom11}, \ref{mom21}, \ref{mom31}, \ref{mom22}), which is obviously simpler. In particular, if one of these relations is not satisfied then the given quadruple of positive numbers $r_i$ cannot be realized as the radii of any Steiner $4$-chain. As an illustration we present a simple example, where our algorithm easily yields a negative answer. 
	
	{\bf Example.}  Given a quadruple of radii $\rho = (1, 2, 3, 4)$, then the bends are $\varkappa = (1, 0.5, 0.33(3), 0.25).$ So $(I_1,I_2,I_3)=(2.083,1.4214,1.1766)$. This implies that (\ref{mom11}) is not satisfied and we conclude that the quadruple $(1, 2, 3, 4)$ is not realizable as the radii of a Steiner $4$-chain.

\end{document}